\documentclass[11pt,draft]{amsart}
\usepackage{mathrsfs}
\usepackage{amssymb}
\usepackage{amsmath}
\usepackage{amsfonts}
\usepackage{amsfonts,enumerate}
\usepackage{pifont}
\usepackage{enumerate}
\usepackage{graphics}
\usepackage{verbatim}
\usepackage{amsthm}

\numberwithin{equation}{section}
\newtheorem{theorem}{Theorem}[section]

\newtheorem{lemma}{Lemma}[section]

\newtheorem{remark}{Remark}[section]

\newcommand{\A}{{\mathcal A}}

\newcommand{\8}{\infty}
\newcommand{\el}{\ell}

\newcommand{\be}{\begin{eqnarray*}}
\newcommand{\ee}{\end{eqnarray*}}
\newcommand{\beq}{\begin{equation}}
\newcommand{\eeq}{\end{equation}}
\newcommand{\beqn}{\begin{equation*}}
\newcommand{\eeqn}{\end{equation*}}
\newcommand{\bs}{\begin{split}}
\newcommand{\es}{\end{split}}

\numberwithin{equation}{section}

\begin{document}

\title{A Littlewood-Paley type theorem for Bergman spaces}


\author{Zeqian Chen}

\address{Wuhan Institute of Physics and Mathematics, Chinese
Academy of Sciences, 30 West District, Xiao-Hong-Shan, Wuhan 430071, China}
\thanks{This work was partially supported by NSFC grant No. 11171338.}

\author{Wei Ouyang}

\address{Wuhan Institute of Physics and Mathematics, Chinese
Academy of Sciences, 30 West District, Xiao-Hong-Shan, Wuhan 430071, China
and Graduate University of Chinese Academy of Sciences, Beijing 100049, China}

\date{}
\maketitle

\markboth{Z. Chen and W. Ouyang}
{Bergman spaces}

\begin{abstract}
In this paper, we prove that the original Littlewood-Paley $g$-functions can be used to characterize Bergman spaces as well.
\end{abstract}


\section{Introduction}\label{intro}

Let $\mathbb{D}$ be the unit disk in the complex plane $\mathbb{C}$ with $\mathbb{T}: = \partial \mathbb{D}$ being the unit circle. Recall that for $0 < p < \8,$ the Hardy space $\mathcal{H}^p$ on $\mathbb{D}$ is defined as the set of all analytic functions $f$ on $\mathbb{D}$ satisfying
\be
\| f \|_{\mathcal{H}^p} : = \sup_{0 \le r < 1} \left ( \int^{2 \pi}_0 | f ( r e^{\mathrm{i} \theta} ) |^p \frac{d \theta}{2 \pi} \right )^{\frac{1}{p}} < \8.
\ee
It is classical that for any $f \in \mathcal{H}^p,$ almost everywhere on $\mathbb{T}$ there exist radial limits $\lim_{r \to 1^-} f (r e^{ \mathrm{i} \theta}),$ denoted by $f (e^{\mathrm{i} \theta}),$ and there holds the relation $\| f \|_{\mathcal{H}^p} = \| f \|_{L^p (\mathbb{T})}.$ We refer to \cite{Zygmund2002} for theory of classical Hardy spaces.

Suppose $f \in \mathcal{H}^p$ and $f = \sum_n a_n z^n$ is the power series of $f.$ Consider the following two quantities
\beq\label{eq:g-seriesfunct}
d (f) (z) = \Big ( \sum^{\8}_{n=0} \big | \triangle_n (f) (z) \big |^2 \Big )^{\frac{1}{2}}
\eeq
where $\triangle_0 (f)(z) = a_0$ and $\triangle_n (f)(z) = \sum_{2^{n-1} \le k < 2^n} a_k z^k$ for $n \ge 1,$ and
\beq\label{eq:g-funct}
g(f) (z) = \left ( \int_0^1  (1-r^2)  | f^{\prime} (r z) |^2 d r \right )^{\frac{1}{2}},
\eeq
for $z \in \mathbb{D} \cup \mathbb{T}.$ Two results on Hardy spaces, essentially due to Littlewood and Paley \cite{LP1936}, assert that
\beq\label{eq:LPTHardy}
\| f \|_{L^p (\mathbb{T})} \approx \| d (f) \|_{L^p(\mathbb{T})} \quad \text{and}\quad \| f - f(0) \|_{L^p (\mathbb{T})} \approx \| g(f) \|_{L^p(\mathbb{T})},
\eeq
the constants involved being only dependent on $p$ with $0< p < \8.$ The two equivalent relations \eqref{eq:LPTHardy} can be considered to be the beginning of the Littlewood-Paley theory.

The main purpose of this paper is to prove that these two equivalent relations \eqref{eq:LPTHardy} hold true as well in the case of $L^p (\mathbb{D})$ replacing $L^p (\mathbb{T}),$ characterizing the so-called Bergman spaces. Recall that for $0 < p < \8,$ the Bergman space $\mathcal{A}^p$ consists of functions $f$ analytic in $\mathbb{D}$ with
\be
\|f\|_{\A^p} =\left ( \int_{\mathbb{D}}|f(z)|^p d A (z) \right )^{\frac{1}{p}}<\infty
\ee
where $d A (z) = d x d y/ \pi$ with $z = x + \mathrm{i}y$ in $\mathbb{D}.$ Note that for $1 \le p < \8,$ $\mathcal{A}^p$ is a Banach space under the norm $\| f \|_{\A^p}.$ If $0 < p <1,$ the space $\mathcal{A}^p$ is a quasi-Banach space with $p$-norm $\| f \|^p_{\A^p}.$

{\it Notion.}\; For two nonnegative (possibly infinite) quantities $X$ and $Y,$ by $X \lesssim Y$ we mean that there
exists a positive constant $C$ depending only on $p$ such that $ X \leq C Y,$ and by $X \thickapprox Y$ that $X \lesssim Y$ and $Y \lesssim X.$

\section{Main results}\label{result}

We state our main results as Theorems \ref{th:LPseries} and \ref{th:LPg}.

\begin{theorem}\label{th:LPseries}
Let $0 < p < \8.$ There are two constants $A_p$ and $B_p$ depending only on $p$ such that
\beq\label{eq:LPseries}
A_p \| f \|_{\A^p} \le \| d(f) \|_{L^p (\mathbb{D})} \le B_p \| f \|_{\A^p}
\eeq
for any $f \in \mathcal{A}^p.$
\end{theorem}

This characterization of those functions in $\A^p$ is a straightforward consequence of the first equivalent relation in \eqref{eq:LPTHardy}, but one of the important features of this characterization is that linear operators obtained by multipliers $m_k$ (of the coefficients $a_k$) that vary boundedly on the dyadic blocks $\triangle_n$ preserve the class $\A^p.$ More generally, this yields a Marcinkiewicz multiplier theorem for Bergman spaces stating that, for any $0 < p < \8$ there exists a constant $C_p$ depending only on $p$ such that
\be
\Big \| \sum_{k=0}^{\8} m_k a_k z^k \Big \|_{\A^p} \le C_p \Big ( \sup_k |m_k| + \sup_{n \ge 0} \sum_{2^n \le k < 2^{n+1}} |m_{k+1} - m_k | \Big ) \| f \|_{\A^p}.
\ee
The proof of this inequality can be obtained as in the case of Hardy spaces (see for example \cite{Zygmund2002}, Theorem XV.4.14).

\

{\it Proof of Theorem \ref{th:LPseries}}.\; Let $0 < p < \8.$ Denote by $f_r ( z) = f ( r z)$ for $0 < r < 1$ and $z \in \mathbb{D}.$ By the first equivalent relation in \eqref{eq:LPTHardy}, one has
\be\begin{split}
\int_{\mathbb{D}} |f(z) |^p d v (z) = & 2 \int_0^1 r  d r \int^{2 \pi}_0 | f_r (e^{\mathrm{i} \theta})|^p \frac{d \theta}{2 \pi} \\
\thickapprox & 2 \int_0^1 r d r \int_0^{2\pi} | d (f_r ) ( e^{\mathrm{i} \theta}) |^p \frac{d\theta}{2\pi}\\
= & \| d (f ) \|^p_{L^p(\mathbb{D})}.
\end{split}\ee
This completes the proof of \eqref{eq:LPseries}.
\hfill$\Box$

\begin{theorem}\label{th:LPg}
Let $0 < p < \8.$ There are two constants $\alpha_p$ and $\beta_p$ depending only on $p$ such that
\beq\label{eq:LPg}
\alpha_p \| f \|_{\A^p} \le \| g(f) \|_{L^p(\mathbb{D})} \le \beta_p \| f \|_{\A^p}
\eeq
for any $f \in \mathcal{A}^p$ with $f(0) = 0.$ Consequently,
\beq\label{eq:LPgpge1}
\| f \|_{\A^p} \approx |f(0)| + \| g(f)\|_{L^p(\mathbb{D})} \quad \text{for}\; 1 \le p < \8,
\eeq
and
\beq\label{eq:LPgp<1}
\| f \|^p_{\A^p} \approx |f(0)|^p + \| g(f)\|^p_{L^p(\mathbb{D})} \quad \text{for}\; 0 < p < 1.
\eeq
\end{theorem}

We will deduce this theorem from some classical results, essentially due to Littlewood and Paley, Marcinkiewicz and Zygmund, and a theorem of Coifman and Rochberg \cite{CR} on atomic decomposition for Bergman spaces (see Lemma \ref{le:complexatomicdecomp} below). The proof is thus considerably elementary.

\begin{lemma}\label{le:complexatomicdecomp} (cf. \cite{Pav2004}, Theorem 8.3.1)\;
Let $0 < p \le 1.$ There exists a sequence $\{a_k\}$ in $\mathbb{D}$ and a constant $C$ such that
$\mathcal{A}^p$ consists exactly of functions of the form
\beq\label{eq:atomdecomp}
f(z)=\sum^{\infty}_{k=1}c_{k}\frac{1-|a_k|^2}{\big ( 1- \bar{a}_{k} z \big )^{2/p +1}},\quad  z\in\mathbb{D}\; ,
\eeq
where $\{c_k\}$ belongs to the sequence space $\el^p$ and the series converges in the quasi-norm topology of $\mathcal{A}^p,$ and
\be
C^{-1} \Big ( \sum_{k}|c_{k}|^p \Big )^{\frac{1}{p}} \le \| f \|_{\A^p} \le C  \Big ( \sum_{k}|c_{k}|^p \Big )^{\frac{1}{p}}.
\ee
\end{lemma}

 {\it Proof of Theorem \ref{th:LPg}}.\; We begin with the first inequality in \eqref{eq:LPg}. Denote by $f_r ( z) = f ( r z)$ for $0 < r < 1$ and $z \in \mathbb{D}.$ Then by the second equivalent relation in \eqref{eq:LPTHardy} we have for any $0 < p < \8,$
\be\begin{split}
\int_{\mathbb{D}} |f(z) - f(0)|^p d v (z)
= & 2 \int_0^1 \| f_r - f_r (0) \|^p_{\mathcal{H}^p} r d r\\
\thickapprox & \int_0^1 r d r \int_0^{2\pi}
\left ( \int_0^1(1 - s^2 )| f_r '(s e^{\mathrm{i} \theta}) |^2 d s \right )^{\frac{p}{2}}\frac{d\theta}{2\pi}\\
\lesssim & \int_0^1 r d r \int_0^{2\pi}
\left ( \int_0^1(1 - s^2 )| f '(r s e^{\mathrm{i} \theta}) |^2 d s \right )^{\frac{p}{2}}\frac{d\theta}{2\pi}\\
\thickapprox & \| g (f ) \|^p_{L^p(\mathbb{D})}.
\end{split}\ee
This proves the first inequality in \eqref{eq:LPg}.

To prove the second inequality in \eqref{eq:LPg} for the case $0 < p \leq 1,$ we will adopt Lemma \ref{le:complexatomicdecomp}. To this end, we write
\be
f_k(z)= \frac{1-|a_k|^2}{\big ( 1- z \bar{a}_{k} \big )^{2/p+1}}.
\ee
An immediate computation yields that
\be\begin{split}
| f_k ' (r z) |^2 =  (2/p +1)^2 |\bar{a}_{k}|^2 (1-|a_k|^2)^2 \frac{1}{ |1 - r z \bar{a}_{k} |^{2(2/p+2)}}
\end{split}\ee
Also, it is easy to check that
\be
|1-t z| \le (1-t) + |1-z| \le 3 |1-t z|,\quad 0 < t \le 1,\; \forall z \in \mathbb{D}.
\ee
Then
\be\begin{split}
g( f_k) (z) = & |a_k | (2/p +1) (1-|a_k|^2) \left( \int_0^1 \frac{(1- r^2) d r}{ |1 - r z \bar{a}_{k} |^{2(2/p +2))}} \right )^{\frac{1}{2}}\\
\lesssim &  (1-|a_k|^2) \left ( \int_0^1 \frac{d r}{\big [ (1-r) + |1 - z \bar{a}_{k} | \big]^{2(2/p +1)+1}} \right )^{\frac{1}{2}}\\
\lesssim &  (1-|a_k|^2) \frac{1}{ | 1- z \bar{a}_{k} |^{2/p +1}}.\\
\end{split}\ee
Hence, for $f = \sum_k c_k f_k$ with $\sum_k | c_k |^p < \8$ we have
\be\begin{split}
\int_{\mathbb{D}} | g (f)(z)|^p d v (z) \leq & \sum_{k=1}^{\infty} |c_k|^p \int_{\mathbb{D}} | g ( f_{k})(z) |^p d v (z)\\
\lesssim & \sum_{k=1}^{\infty} |c_k|^p(1-|a_k|^2)^p\int_{\mathbb{D}} \frac{1}{| 1 - z \bar{a}_{k} |^{2 + p}} d v(z)\\
\lesssim & \sum_{k=1}^{\infty}|c_k|^p,
\end{split}\ee
where the last inequality is obtained by the fact that
\be
\int_{\mathbb{D}} \frac{1}{| 1 - z \bar{w}|^{2 +  p}} d v (z) \thickapprox \frac{1}{(1 - |w|^2)^p} \quad \text{as}\; |w| \to 1^-,
\ee
for $p > 0$ (see Theorem 1.7 in \cite{HKZ}). By Lemma \ref{le:complexatomicdecomp}, we conclude the second inequality in \eqref{eq:LPg} for the case $0 < p \le 1.$

Finally, let $1 < p < \8.$ If $f \in \mathcal{A}^p,$ then $f$ has the integral representation
\be
f(z) = \int_{{\mathbb{D}}} \frac{f(w) d v (w)}{(1- z \bar{w} )^2},\quad \forall z \in \mathbb{D}.
\ee
An immediate computation yields that
\be
\left | f' (r z) \right | \lesssim (1-|r z|^2)^{ - \frac{1}{2}} \int_{\mathbb{D}} \frac{|f(w)| d v (w)}{|1-r z \bar{w} |^{\frac{5}{2}}}.
\ee
Then,
\be\begin{split}
g( f)^2 (z) \lesssim & \int_0^1 \frac{1-r^2}{1-|r z|^2} \left | \int_{{\mathbb{D}}} \frac{|f(w)| d v (w)}{|1 - r z \bar{w} |^{\frac{5}{2}}} \right |^2 d r \\
\lesssim & \int_{\mathbb{D} \times \mathbb{D}}|f(w) f(\xi) | d v (w) d v (\xi) \int_0^1 \frac{d r}{|1- r z \bar{w} |^{\frac{5}{2}}|1 - r z \bar{\xi} |^{\frac{5}{2}}}\\
\lesssim & \int_{\mathbb{D} \times \mathbb{D}} |f(w) f(\xi)| d v (w) d v (\xi)\\
 & \quad \times \int_0^1 \frac{d r}{\big [ |1 - z \bar{w} | + (1-r) \big ]^{\frac{5}{2}} \big [ |1- z \bar{\xi} | + (1-r) \big ]^{\frac{5}{2}}}\\
\lesssim & \int_{\mathbb{D} \times \mathbb{D}} |f(w) f(\xi) | d v(w) d v(\xi)\\
 & \quad \times \left ( \int_0^1 \frac{d r}{\big [ |1- z \bar{w} | + (1-r) \big ]^5} \int_0^1 \frac{d r}{\big [ |1- z \bar{\xi}  | + (1-r) \big ]^5} \right)^{\frac{1}{2}}\\
\lesssim & \int_{\mathbb{D} \times \mathbb{D}} \frac{|f(w) f(\xi) |}{|1- z \bar{w} |^2 | 1 - z \bar{\xi} |^2} d v (w) d v (\xi)\\
= & \left ( \int_{\mathbb{D}} \frac{|f(w)|}{|1- z \bar{w} |^2}d v (w) \right )^2.
\end{split}\ee
However, the mapping
\be
f \longmapsto \int_{\mathbb{D}} \frac{f(w)}{|1 - z \bar{w} |^2} d v (w)
\ee
is bounded on $L^p (\mathbb{D})$ for $1< p < \8$ (e.g. Theorem 1.9 in \cite{HKZ}). Therefore, we conclude the second inequality in \eqref{eq:LPg} for the case $1 < p < \8.$
\hfill$\Box$

\begin{remark}\label{rk:gqfunct}
\begin{enumerate}[{\rm (1)}]

\item Since 1930's the Littlewood-Paley theory was developed considerably and mainly carried out by E. M. Stein \cite{Stein1970}, widening its applicability both in the
classical setting involving $\mathbb{R}^n$ (even when $n = 1$) and in abstract situations involving, among other things, Lie groups, symmetric spaces, diffusion
semigroups and martingales. We consult \cite{FJW} and references therein for more recent information.

\item Some real-variable characterizations of Bergman spaces involving maximal and area integral functions in terms of the Bergman metric, have been obtained recently by the present authors \cite{CO1, CO2}.

\end{enumerate}
\end{remark}


\begin{thebibliography}{99}

\bibitem{CO1}Z. Chen and W. Ouyang,
Maximal and area integral characterizations of Bergman spaces in the unit ball of $\mathbb{C}^n,$
arXiv: 1005.2936.

\bibitem{CO2}Z. Chen and W. Ouyang,
Real-variable characterizations of Bergman spaces in the unit ball of $\mathbb{C}^n,$
arXiv: 1103.6122.

\bibitem{CR}R. Coifman and R. Rochberg,
Representation theorems for holomorphic and harmonic functions in $L^p,$
{\it Asterisque} {\bf 77} (1980), 11-66.

\bibitem{FJW}M. Frazier, B. Jawerth, and G. Weiss,
{\it Littlewood-Paley Theory and the Study of Function Spaces,}
American Mathematical Society, 1991.

\bibitem{HKZ}H. Hedenmalm, B. Korenblum, and K. Zhu,
{\it Theory of Bergman Spaces,}
Springer-Verlag, New York, 2000.

\bibitem{LP1936}J. E. Littlewood and R. E. A. C. Paley,
Theorems on Fourier series and power series II,
{\it Proc. London Math. Soc.} {\bf 42} (1936), 52-89.

\bibitem{Pav2004}M. Pavlovi\'{c},
{\it Introduction to Function Spaces on the Disk,}
Matemati\v{c}cki institut SANU, Beograd, 2004.

\bibitem{Stein1970}E. M. Stein,
{\it Topics in harmonic analysis related to the Littlewood-Paley theory,}
Princeton University Press, Princeton, New Jersey, 1970.

\bibitem{Zygmund2002}A. Zygmund,
{\it Trigonometric Series (Third Edition),}
Cambridge University Press, Cambridge, England, 2002.

\end{thebibliography}
\end{document}